\documentclass[12pt]{article}
\usepackage{tikz}

\usepackage{graphicx}
\usepackage{color}
\usepackage{amssymb}
\usepackage{amsmath}

\newtheorem{theorem}{Theorem}

\newtheorem{lemma}[theorem]{Lemma}

\usepackage{graphicx}
\usepackage{color}

\makeatletter
\@addtoreset{equation}{section}
\makeatother

\newcommand{\mf}[1]{{\mathfrak #1}}

\newcommand{\reff}[1]{(\ref{#1})}

\newcommand{\Cal}[1]{{\mathcal #1}}

\newcommand{\bL}{\bold L}
\newcommand{\bA}{\bold A}
\newcommand{\bY}{\bold Y}
\newcommand{\bV}{\bold V}
\newcommand{\bW}{\bold W}

\def\p1ff{<\mf f, \mf L_s \mf f>}
\def\p1gg{<\mf g, \mf L_s \mf g>}

\def\emptysquare{{\hbox{\vrule height6pt width0.6pt depth0pt%
\vbox{\hrule height0.6pt width4.8pt depth0pt%
\vglue4.8pt%
\hrule height0.6pt width4.8pt depth0pt}%
\vrule height6pt width0.6pt depth0pt}}}

\def\qed{\unskip\nobreak
\hfil\penalty50\hskip1.75em\null\nobreak\hfil\emptysquare
{\parfillskip=0pt \finalhyphendemerits=0 \par}\medskip}

\newenvironment{demo}{\noindent {\bf Proof:}~}{\qed \medskip}
\begin{document}

\title{Some considerations on the back door theorem and conditional randomization.}
\author{ {Julieta Molina$^{a,b}$ }, { Lucio Pantazis$^{a}$ } , {Mariela Sued$^{a,b}$}\footnote{Corresponding author. Email: msued@dm.uba.ar. Address: Intendente Guiraldes 2160 - Ciudad Universitaria - C1428EGA -  Telephone and fax number (++54 +11) 4576-3375}\\
{\small $^a$ \sl  Universidad de Buenos Aires, \small $^b$ \sl CONICET, Argentina} }
\date{~}
\maketitle


\begin{abstract}
In this work we   propose a different
 \textit{surgical modified model} for  the construction of counterfactual variables under non parametric structural equation models. This  approach allows  the simultaneous representation of counterfactual responses and observed treatment assignment, at least when the intervention is done in one node. Using the new proposal, the $d$-separation criterion is used verify conditions related with  ignorability or conditional ignorability and a new proof of the back door theorem is provided under this framework.

\end{abstract}

\section{Introduction}

 In the potential outcome framework (Rubin, 1974),
   identifiability of the average treatment effect is guaranteed
 under the assumption of  exchangeability (or ignorability). Let  $Y_t$ and
 $Y_c$ denote the potential outcomes  under treatment  level $t$ and
 $c$, respectively. More precisely, $Y_t$ is the outcome variable that would have been observed in a hypothetical world in which all individuals
 received treatment level $t$, and $Y_c$ is the outcome variable that would have been observed if all individuals were treated under level $c$.
 If $A$ denotes the observed treatment assignment,under exchangeability potential
 outcomes are independent of treatment assignment: $Y_a \coprod A$, for $ a=t,c$.
If $Y$ denotes the observed response, this assumption together with
consistency ($Y=Y_A$) and positivity ($0<P(A=t)<1$), guarantees that
$E[Y_a]=E\left[Y\vert A=a\right]\;,\text{for $a=c,t$}$
and so, the average treatment effect is identified by the
distribution of observed data $(A,Y)$, by the formula
$\text{ATE}=E[Y_t]-E[Y_c]= E\left[Y\vert  A=t\right]-E\left[Y\vert
A=c\right]$.

When exchangeability is not a reasonable assumption,
with the help of an observed vector $L$, conditional exchangeability still
guarantees identifiability. Conditional exchangeability states that
potential  outcomes are independent of treatment assignment,
given   $L$: $ Y_a \coprod A\vert L$, for $a=c,t$.

This assumption together with consistency and positivity
($P(A=a\vert L=\ell)>0$ if $P(L=\ell)>0)$), ensures that
$
E[Y_a]=E\left[E[Y\vert A=a,L]\right]\;,\text{for $a=c,t$}$,
and so, the average treatment effect is identified by the
distribution of observed data $(L,A,Y)$, by the formula
$\text{ATE}=E[Y_t]-E[Y_c]= E\left[E[Y\vert
A=t,L]\right]-E\left[E[Y\vert A=c,L]\right]$.

On the other hand, the $d$-separation criterion is a graphical tool designed to check
independence and conditional independence between coordinates (or
sub-vectors) of a random vector whose distribution  satisfies the Markov
factorization relative to a given directed acyclic graph (DAG).
Throughout this work, we assume that the reader is familiarized with graphical models, DAG's, and the results relating $d$-separation with conditional independences.  Pearl's book (Pearl, 2009a) contains all background information required to read this work.

Working with non parametric structural equation models (NPSEM),  potential outcomes are defined replacing the equations related to the treatment nodes by  the constants corresponding to the desired  intervention: following Pearl's work (Pearl (2009a)), if we are interested in $Y_t$, we  use the modified model
$M_t$, where the function associated to the node $\bA$  is fixed at $t$: $f_A= t$.

This construction does not allow $A$ (the observed treatment assignment) and $Y_t$
to be represented with the same set of structural equations.
In this work we propose a new approach to construct the counterfactual variables using NPSEM in such a way that the potential outcomes and treatment variables are jointly represented by the same vector, at least when the intervention is done in one node.
 This allows  the use of the $d$-separation rules
to guarantee the identifiability of the distribution of counterfactual variables.   Moreover, when the intervention addresses only one node, the assumptions of the back-door theorem imply conditional randomization.

The twin DAG's presented by Balke and Pearl (1994) also allows to construct simultaneously observed  and counterfactual variables. Richardson and Robins (2013) present a deep discussion about these models. However,  considering that the back door theorem is one of the most popular criterion for identifying the distribution of counterfactual variables, we do care about constructing counterfactual variables and treatment assignments using the graph involved in the back door theorem, namely in the DAG where arrows emerging from nodes associated to the intervention are removed.

We started thinking about DAG's and counterfactual variables  a couple of years ago, inspired in  the back door theorem.
We felt  that  the d-separation criterion could be used to deduce  randomization or conditional randomization, reconciling  different approaches that give rise to the same formula for the distribution of the potential outcomes. We succeeded, at least,  when the intervention is done in only one node and this case has been included in a Final Degree Thesis (Licenciatura en Matematica) (Cacheiro, 2011).

During the last years other authors have been concerned  with this problem. In a working paper,  Richardson and Robins (2013) present a graphical theory based on single world intervention graphs (SWIGs),  unifying causal directed graphs and potential outcomes.

 Although their graphs differ from the graphs in this paper, our perspectives are rather similar. One difference is that although  they assume an underlying  NPSEM ,they do not assume
independence of the disturbances (as we do) and prefer to work with the Finest Fully Randomized Causally Interpretable Structured Tree Graphs  models (FFRCISTG).

In this work, we start assuming  NPSEM with independent errors (NPSEM-IE). Although  NPSEM-IE are less general than FFRCISTG models,  we decided to present our proposal under this setting, considering   that these models are still adopted by a large portion of the community (like Pearl's followers).
However, in Section \ref{FF} we show that all the results presented in this work remain true if FFRCISTG models are assumed instead.


This work is organized as follows. In Section \ref{toyexample} we present a simple example with a three node DAG, explaining the main idea to construct jointly  $A$ and $Y_a$. We check in this example that the relation with the back door theorem and conditional randomization holds. In Section \ref{constant} we generalize these results, first when the intervention is done in just one node and then generalized for many nodes. In Section \ref{FF} we discuss our results under FFRCISTG models.

 To end this introduction we would like to  establish a subtle difference oftenly omitted. Given   a DAG $G$, we use $\bV=\{\bV_1,\bV_2,\cdots,\bV_n\}$ to denote the nodes of a  graph,
 while random variables associated with a given node $\bV_i$ are denoted by $V_i$ or some perturbation of $V_i$, like $V_{i,t}$ or $V_i^t$, as it will be explained later on.

\section{Toy Example - Main Idea}
\label{toyexample}
Assume that the causal diagram  associated with the  problem of interest is given by
the  DAG $G$ \begin{center}
\begin{tikzpicture}[inner sep=1pt,scale=1.5]
  \tikzstyle{every node}=[draw,shape=circle];
  \node (a) at (5,1)     {$\bA$};
  \node (y) at (7,1){$\bY$};
  \node (l) at (6,2)  {$\bL$};
  \draw[->] (l)--(a);
   \draw[->] (a)--(y);
   \draw[->] (l)--(y);
\end{tikzpicture}
\newline
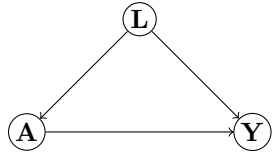
Figure 1: the original DAG $G$.
\end{center}

In terms of NPSEM's (Pearl,  2009b) this means that there exists a set of functions $
F=\{f_\bL, f_\bA,f_\bY\}$ and jointly independent disturbances
$U=\{U_\bL,U_\bA,U_\bY\}$, that give rise to factual  variables
according to the following recursive system
\begin{equation}
\label{model}
L=f_\bL(U_\bL)\;,\quad A=f_\bA(L,U_\bA)\;,\quad Y=f_\bY(L,A,U_\bY)\;.
\end{equation}
We use $M=(  F, U)$ to denote the model which defines  factual
variables.  To emulate the intervention $do(a)$, Pearl (2009b) considers a
model $M_{a}$ where the function $f_\bA$ is replaced by the constant
$a$, while the disturbances remain unchanged: $M_{a}=(  F_a,U)$,
with $ {F}_a=\{f_{a,\bold L},f_{ a,\bold A},f_{a,\bold Y}\}$, where
 \begin{equation}
 f_{a,\bL}=f_\bL\;,\quad f_{ a,\bA}=a\;,\quad f_{a,\bY}=f_\bY
\;.
\end{equation}

 Variables obtained iterating the functions in model $M_{a}$
 using the same vector of disturbances $U=\{U_\bL,U_\bA,U_\bY\}$ are
denoted with subindex $a$: $L_a$, $A_a$ and $Y_a$. In this way, the
counterfactual response of interest at level $a$ is given by
$Y_{a}$.

Our proposal  to represent counterfactual variables consist in the use of
a new system of functions, in which  the value $a$ is inserted in lieu of
the variable corresponding to the node
every time this one is required  by the
recursion. To do so we change the functions related to
each node that have
as parent. In the present
example, $M^{a}=(  F^a,U)$, with $
 F^a=\{f^{a}_\bL,
f^{a}_\bA,f^{a}_\bY\}$, where
\begin{equation}
f^{a}_\bL=f_\bL\,\quad  f^{a}_\bA=f_\bA\;,\quad
f^{a}_\bY(\ell,u)=f_\bY(\ell,a,u)\;.
\end{equation}
Note that this new set of functions is compatible with the DAG $G_{\underline \bA}$, where arrows emerging from $\bA$ are removed:
 \begin{center}
\begin{tikzpicture}[inner sep=1pt,scale=1.5]
  \tikzstyle{every node}=[draw,shape=circle];
  \node (x) at (0,1)     {$\bA$};
  \node (y) at (2,1){$\bY$};
  \node (s) at (1,2)  {$\bL$};
    \draw[->] (s)--(x);
   \draw[->] (s)--(y);
\end{tikzpicture}
\newline
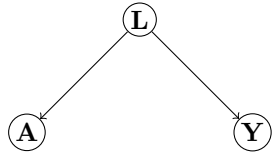
Figure 2: $G_{\underline \bA}$ , constructed removing in  $G$ arrows emerging from $\bA$.
\end{center}

Variables constructed iterating the functions in $F^a$ and using the same vector of disturbances $\{U_\bL,U_\bA,U_\bY\}$ are denoted with supraindex $a$: $L^a$, $A^a$ and $Y^a$. Then, we get that the distribution of $(L^a, A^a, Y^a)$ is compatible with $G_{\underline \bA}$.

The following Lemma summarized the main results of this section.
\begin{lemma}
\begin{enumerate}
\item \label{uno}$L_a=L^{a}=L$,  $A^{a}=A$ and $Y_{a}=Y^{a}$.
\item \label{dos} $\bA$ and $\bY$ are $d$-separated by $\bL$ in  $G_{\underline \bA}$ and so, since
the distribution of $(L^a,A^a,Y^a)$ is compatible with $G_{\underline \bA}$, we get that $A^a$ is independent of $Y^a$ given $L^a$.
\item From the previous results,  we conclude that $A$ is independent of $Y_a$ given $L$. Thus,
conditional randomization holds.
    \end{enumerate}
    \end{lemma}
\section{Intervention with constant regimes}
\label{constant}
\subsection{First step: one node }


Consider a causal DAG $G$ with nodes $\bV_1,\cdots,\bV_n$, labeled
in a compatible way with $G$. Recall that in the graph terminology, we say that
$\bV_i$ is a parent of $\bV_j$ if an arrow points from $\bV_i$ to $\bV_j$. We use $PA_G(\bV_j)$ to denote the set of parents of $\bV_j$ in $G$. If  $\bV_i$ has a directed path to $\bV_k$ we say that  $\bV_i$ is an ancestor of $\bV_k$, and use $An_G(\bV)$ to denote the set of ancestors of $\bV$ in $G$.

Consider a collection
of independent random variables $U=\{U_1,\cdots, U_n\}$. Let $\Cal
V_i$ denote  the common support of any random variables associated
with the node $\bV_i$ and let   $\Cal U_i$ denote the support of
$U_i$. A set of functions  $F=\{f_i:i\geq 1\}$ is said to be  compatible with $G$
if  for each $i=1,\cdots n$ we get that
\begin{equation}
\label{fi} f_i: \prod_{\bV_j \in  PA_G(\bV_i)}\Cal  V_j \quad\times\quad \Cal  U_i \quad\longrightarrow \quad
\Cal  V_i\;.
\end{equation}

Given a set  $F=\{f_i:i\geq 1\}$ of compatible functions
with $G$, and independent $U=\{U_1,\cdots, U_n\}$,  factual
variables are defined by the recurrence
$$V_i=f_i(PA_i,U_i)\;,$$
where $PA_i$ are the random variables (already defined by the
recurrence) associated with the nodes in $PA_G(\bV_i)$.
Note that, by construction, the distribution of $(V_1,\dots,V_n)$ is compatible with $G$, meaning that it satisfies the Markovian factorization induced by $G$.
We use
$M=(F,U)$ to denote the model that gives rise to factual variables.

In order to represent an intervention at level $a$ for a given node
$\bA$, Pearl (2009b) defined the "Surgically modified model" $M_a=(F_a,U)$, considering $F_a=\{f_{a,i}:i\geq 1\}$, where $f_{a,i}=f_i$ if
$\bV_i\not=\bA$ and for $\bV_j=\bA$, $f_{a,j}= a$. Counterfactual
variables are defined by this new set of functions and the same
disturbances $\{U_1,\cdots, U_n\}$, by the recurrence
$$V_{a,i}=f_{a,i}(PA_{a,i},U_i)\;,$$
where $PA_{a,i}$ are the random variables (already defined by the
recurrence) associated with the nodes in $PA_G(\bV_i)$.

Before presenting our proposal for constructing counterfactual variables, recall that given a DAG $G$ and a node $\bA$ in $G$, $G_{\underline \bA}$ is the graph obtained by removing from $G$  all arrows emerging from $\bA$.  We will now introduce a new set of functions  $F^a=\{f^a_{i}:i\geq 1\}$, compatible with $G_{\underline \bA}$,  that will allow the simultaneous definition of
both the observed assignment random variable $A$ associated with
the node $\bA$,  and the counterfactual responses.
To achieve this, if $\bA\notin PA_G(\bV_i)$ we get that $PA_{G_{\underline \bA}}(\bV_i)=PA_{G}(\bV_i)$  and define   $f^a_{i}$ being equal to $f_i$ . When $\bA\in PA_G(\bV_i)$, $f^a_{i}$ is obtained fixing the value $a$ at the original function $f_i$.  To be more precise,  if $\bA\in PA_G(\bV_i)$, without loss of generality, we can assume that
\begin{equation}
\label{fioriginal} f_i: \prod_{\bV_j \in  PA_G(\bV_i)\setminus\bA}\Cal  V_j \quad\times\quad \Cal A\quad\times\quad  \Cal  U_i \quad\longrightarrow
\quad\Cal  V_i\;,
\end{equation}
where $\Cal A$ denotes the set of possible values to be taken by variables associated to node $\bA$. Since   $PA_{G_{\underline \bA}}(\bV_i)=PA_{G
}(\bV_i)\setminus \bA$ and  $F^a$ should be  compatible with  $G_{\underline \bA}$,  we need $f_i^a$ to satisfy the following condition:
\begin{equation}
\label{fi} f^a_i: \prod_{\bV_j \in  PA_G(\bV_i)\setminus\bA}\Cal  V_j \;\quad\times\quad  \Cal  U_i \quad\longrightarrow
\quad\Cal  V_i\;.
\end{equation}

Then,  for $\overline v_i\in \prod_{\bV_j \in  PA_G(\bV_i)\setminus\bA}\Cal  V_j $, $u\in \Cal U_i$,  we define
$$f^a_i(\overline v_i, u)=f_i(\overline v_i,a, u).$$

Let $V_i^a$ denote the variables obtained by the
recurrence based on these new functions:
$$V^a_{i}=f^a_{i}(PA^a_{i},U_i)\;,$$
where $PA^a_{i}$ are the random variables (already defined by the
recurrence) associated with the nodes in $PA_{G_{\underline \bA}}(\bV_i)$. Note that the distribution of $(V_1^a,\dots, V_1^n)$ is compatible with
$G_{\underline \bA}$. Let
$M^a=(F^a,U)$.

The following Lemma explains how variables defined under models  $M$, $M_a$ and $M^a$ are related.
\begin{lemma}
\label{2}  The random variables
associated with both modified models $M_a=(F_a,U)$ and $M^a=(F^a,U)$
are the same, with the exception of those associated with node $\bA$:
$${V_{i,a}=V^a_i\quad \text{if \quad $\bV_i\not=\bA$}}\;.$$
Variables associated with the node $\bA$ defined by $M=(F,U)$ and
$M^a=(F^a,U)$, respectively, are equal:
$$A=A^a\;.$$
Moreover, if $\bV_i$ is not a descendent of $\bA$, we get that
$$ V_i=V_{a,i}=V_i^a\;.$$
\end{lemma}

To end this section, we state the back door Theorem (Pearl, 2009b), and provide a new proof of it.
\begin{theorem}
\label{bd}\textbf{The Back Door Criterion}
Consider a set of nodes $\bL\subset \{\bV_1,\cdots, \bV_n\}$, {such that  $\,\,\,\bL\cap \bA=\emptyset$}. Assume that the following conditions hold:

\begin{enumerate}
\item \label{desc}No element of $\bL
$ is a descendent of $\bA$ in $G$,
\item \label{block} $\bL$  blocks all back door paths from $\bA$ to $\bY$ in $G$.
\end{enumerate}
Then, $Y_a$ in independent of $A$ given $L$ and so $$P(Y_a=y)=\sum _\ell P(Y=y \vert A=a,L=\ell)P(L=\ell)\;.$$
\end{theorem}
\bigskip
\begin{demo}
To prove that conditional ignorability holds, meaning that $Y_a$ in independent of $A$ given $L$, we note that under the assumption of Theorem \ref{bd},  considering the  results presented in Lemma \ref{2}, we get that
\begin{enumerate}\item \label{desc} If no element of $\bL
$ is a descendent of $\bA$ in $G$, then $L=L_a=L^a$.

\item \label{block}If  $\bL$  blocks all back door paths from $\bA$ to $\bY$ in $G$, then $\bA$ and $\bY$ are $d$-separated by $\bL$ in
$G_{\underline \bA}$, and so $A^a$ and $Y^a$ are independent given $L^a$.
\end{enumerate}
Finally, applying again to the results stated in  Lemma \ref{2}, we also know that $A^a=A$ and $Y^a=Y_a$. So, if $\bL$ satisfies both conditions \reff{desc}
and \reff{block}, we can conclude that  $Y_a$ is
independent of $A$ given $L$. This means that conditional ignorability holds, as we wanted to prove. Thus, under positivity, the distribution of the counterfactual variables can be identified:
$$P(Y_a=y)=\sum _\ell P(Y=y \vert A=a,L=\ell)P(L=\ell)\;.$$
\end{demo}

\subsection{Intervention with constant regimes - many nodes }
\label{constantemuchos}

Assume now that we wish to intervene in a set of nodes
$\bA_{set}=\{\bA_1,\dots,\bA_k\}$. Consider  $a_i\in \Cal A_i$,
where $\Cal A_i$ denotes the support of variables associated with node
$\bA_i$, and  let $a=(a_1,\cdots,a_k)$. Following the new surgically
modified model, we will change the functions  related to those nodes
whose parents have some $\bA_j$.

As in the one node case, given a DAG $G$, let $M=(F, U)$ denote the model (compatible with $G$) for observed variables $(V_1,\cdots,V_n)$.
Let $(V_{a,1},\dots, V_{a,n})$ denote the vector of variables determined by the model $M_a=(F_a,U)$ proposed by Pearl, with  $F_a=\{f_{a,i}:i\geq 1\}$, where $f_{a,i}=f_i$ if
$\bV_i$ does not belong so the set $\bA_{set}$,  and when $\bV_j= \bA_i$ for some $i$, $f_{a,j}= a_i$.

We will now generalize our  construction
presented for single node intervention in this new scenario. To do so, we consider $M^a=(F^a,U)$, for $F^a=\{f^a_{i}:i\geq 1\}$ compatible with $G_{\underline {\bA_{set}}}$, the graph obtained   removing in $G$ all arrows emerging from the set $\bA_{set}$.
Note that the set of parents of a given node $\bV_i$ in $G_{\underline {\bA_{set}}}$ is obtained eliminating from the set of parents  of $\bV_i$ in the original DAG $G$ all nodes in $\bA^i_{set}=\bA_{set}\cap PA_G(\bV_i)$, namely, we have that $PA_{G_{\underline {\bA_{set}}}}(\bV_i)=PA_G(\bV_i)\setminus \bA^i_{set}$. Therefore, the definition of  $f_i^a$ depends on whether the set  $\bA_{set}^i=\bA_{set}\cap PA_G(\bV_i)$ is empty or not.  Now, if $\bA_{set}^i=\emptyset$, we get that $PA_{G_{\underline {\bA_{set}}}}(\bV_i)=PA_G(\bV_i)$ and we define $f_i^a=f_i$. When
 $\bA_{set}^i\not=\emptyset$, we can assume that
  \begin{equation}
\label{fioriginalmultiple} f_i: \prod_{\bV_j \in  PA_G(\bV_i)\setminus\bA_{set}^i}\Cal  V_j\quad \times\quad \prod_{\bA_j \in  \bA_{set}^i} \Cal A_j\quad \times\quad  \Cal  U_i \quad\quad\longrightarrow\quad\quad
\Cal  V_i\;,
\end{equation}
and consider
     \begin{equation}
\label{fiamultiple} f^a_i: \prod_{\bV_j \in  PA_G(\bV_i)\setminus\bA_{set}^i}\Cal  V_j\quad \times\quad  \Cal  U_i \quad\quad\longrightarrow\quad\quad
\Cal  V_i\;,
\end{equation}
where for $\overline v_i\in \prod_{\bV_j \in  PA_G(\bV_i)\setminus\bA_{set}^i } \Cal {V} _j$, $u\in \Cal U_i$,  we define
$$f^a_i(\overline v_i, u)=f_i(\overline v_i,a^i, u),$$
 putting in $a^i$ all the coordinates of the vector $a=(a_1,\cdots, a_k)$
 corresponding to the set $\bA^i_{set}$: $a^i=(a_j:\bA_j\in \bA^i_{set})$.
In other words, when $PA_G(\bV_i)\cap \bA_{set}\not=\emptyset$ we construct the function $f_i^a$  fixing at $f_i$  the value $a_j$, each time the
value of the variable related to the node $\bA_j$ is required by the original function $f_i$ (for each $j$ such that $\bA_j\in \bA_{set}^i$).

Let $(V_1^a,\dots,V_n^a)$ denote the vector of  variables obtained by the
recurrence based on these new functions  ($F^a$) and disturbances $U$. Once more, we get that the distribution of $(V_1^a,\dots, V_n^a)$ is compatible with $G_{\underline {\bA_{set}}}$.
The results are presented in what follows.

\begin{lemma}
\label{l3}
 \label{mismo}Let $A=(A_1,\cdots,A_k)$ and  $A^a=(A^a_1,\cdots,A^a_k)$ denote the random variables related to the nodes $\bA_1,\cdots,\bA_k$, according to model $M$ and $M^a$, respectively.
{If $\,\,\bW\cap \bA_{set}=\emptyset$}, then the following version of the consistency assumption holds:
$$\{A^a=a\;,\quad \;W^a=w\}=\{A=a\;,\quad W=w\}\;.$$

The random variables associated with both modified models
$M_a=(F_a,U)$ and $M^a=(F^a,U)$ are the same, with the exception of those associated with
nodes in  $\bA_{set}$:
$$V_{i,a}=V^a_i\quad {\text{if $\,\,\bV_i\not\in \bA_{set}$}}\;.$$
\end{lemma}

Finally, we include a new proof of the Back Door Theorem, using the independences deduced from its  assumptions and Lemma \ref{l3}.

\begin{theorem}\textbf{Back Door Criterion: Many Nodes}
Consider a set of nodes $\bL\subset \{\bV_1,\cdots, \bV_n\}$, such that  $\bL\cap \bA_{set}=\emptyset$.  Assume that the following conditions hold:

\begin{enumerate}
\item No element of $\bL$ is a descendent of {$\bA_{set},$}

\item $\bL$ blocks all back door paths from $\bA_{set}$ to $\bY$ in $G$.
\end{enumerate}
Then, $$P(Y_a=y)=\sum _\ell P(Y=y \vert A=a,L=\ell)P(L=\ell)\;,$$
\end{theorem}
with $a=(a_1,..a_k)$.

\begin{demo}
Under the present assumptions we get that

\begin{enumerate}
\item If no element of $\bL$ is a descendent of $\bA_{set}$, {then $L=L_a=L^a.$}

\item If $\bL$ blocks all back door paths from $\bA_{set}$ to $\bY$ in $G$, then $\bA_{set}$ and $\bY$ are $d$-separated by $\bL$ in
$G_{\underline {\bA_{set}}}$, and so $A^a$ and $Y^a$ are independent given $L^a$.

\end{enumerate}

Finally, if $\bL$ satisfies the previous conditions, by Lemma \ref{mismo},
 we get that $\{A^a=a,\;L^a=\ell\}=\{A=a,\;L=\ell\}$ for any $\ell$, $\{Y^a=y,A^a=a,\;L^a=\ell\}=\{Y=y,A=a,\;L=\ell\}$  for any $(\ell, y)$ and
 so, under positivity

\begin{eqnarray*}&&P(Y_a=y)=P(Y^a=y)=\sum _\ell P(Y^a=y \vert L^a=\ell)P(L^a=\ell)=\\
&&
\sum _\ell P(Y^a=y \vert L^a=\ell, A^a=a)P(L^a=\ell)=\sum _\ell P(Y=y \vert A=a,L=\ell)P(L=\ell)\;.
\end{eqnarray*}
\end{demo}

\section{FFRCISTG models}
\label{FF}

In the previous results we have used the rules of $d$-separation to detect independence or conditional independence between variables of a random vector. To do so, given a graph $G$, all we required from the joint distribution of our vector was compatibility with $G$. When variables are constructed following a NPSEM-IE, the Markov  factorization induced by $G$ holds automatically, and that is why our results are valid when the errors are independent.

However, the Markov factorization remains true  under weaker conditions. For instance, let $v=(v_1,\dots ,v_n)\in \prod_{j=1}^n \Cal V_j$ and call  $v_{pa_{G}(\bV_i)}$ the  subvector of $v$ containing the coordinates related with nodes in  the set $PA_G(\bV_i)$, namely $v_{pa_{G}(\bV_i)}=(v_j: \bV_j \in PA_G(\bV_i))$. If
\begin{equation}
\label{FFRCISTG}
\left\{f_i\left(v_{pa_{G}(\bV_i)},U_i\right): \bV_i\in G\right\} \text{ are independent,}\quad \hbox{for all $v\in \prod_{j=1}^n \Cal V_j$}\;,
\end{equation}
then, the distribution of the vector whose variables are constructed  with $M=(F,U)$ is compatible with the graph $G$.
This condition, mainly defines the FFRCISTG models (Richardson and Robins 2013).

It is worthy  to note that if  $M=(F,U)$  satisfies condition (\ref{FFRCISTG}) relative to $G$,  the  intervened model $M^a=(F^a,U)$  defined in Section \ref{constantemuchos}, also satisfies condition (\ref{FFRCISTG}) relative to $G_{\underline \bA_{set}}$, since
$$
\left\{f_i^a\left(v_{pa_{G_{\underline \bA}}(\bV_i)},U_i\right): \bV_i\in G_{\underline \bA}\right\}=
\left\{f_i\left(v^a_{pa_{G}(\bV_i)},U_i\right): \bV_i\in G\right\}$$
where  $v^a_{pa_{G}(\bV_i)}$ denotes the vector that results from replacing $v_j$ with $a_j$ for $\{j:\bA_j \in \bA_{set}\}$. Then,  variables constructed by the model $M^a$ satisfy the factorization induced by the graph $G_{\underline \bA_{set}}$, allowing the use of $d$-separation rules, and thus, extending our results for this new model.

\end{document}